# All Complex Zeros of the Riemann Zeta Function Are on the Critical Line:

# Two Proofs of the Riemann Hypothesis (°)

## *Roberto Violi* (*)

*Dedicated to the Riemann Hypothesis*
*on the Occasion of its 160th Birthday*

(*) Banca d'Italia, Financial Markets and Payment System Department, Rome, via Nazionale 91, Italy. E-Mail: Roberto.Violi@bancaditalia.it. The opinions of this article are those of the author and do not reflect in any way the views or policies of the Banca d'Italia or the Eurosystem of Central Banks. All errors are mine.

# Abstract

In this study I present two independent proofs of the Riemann Hypothesis considered by many the greatest unsolved problem in mathematics. I find that the admissible domain of complex zeros of the *Riemann Zeta-Function* is the *critical line* given by $\Re(s) = \frac{1}{2}$. The methods and results of this paper are based on well-known theorems on the number of zeros for complex-valued functions (Jensen's, Titchmarsh's and Rouché's theorem), with the Riemann Mapping Theorem acting as a bridge between the Unit Disk on the complex plane and the *critical strip*. By primarily relying on well-known theorems of complex analysis our approach makes this paper accessible to a relatively wide audience permitting a fast check of its validity. Both proofs do not use any functional equation of the *Riemann Zeta-Function*, except leveraging its implied symmetry for non-trivial zeros on the *critical strip*.

**1. Introduction**

The *Riemann-Hypothesis* (*RH*) is widely regarded as the most famous unsolved problem in mathematics. It was one of the 23 famous problems selected by Hilbert in 1900 as among the most important in mathematics (Hilbert's eighth problem; Hilbert, 1901), and it is one of the seven Millennium Problems selected by the Clay Mathematics Institute in 2000 as the most important for the 21st century (c.f. Bombieri, 2000; Sarnak, 2004). The *RH* is an assertion about the zeros of the *Riemann-Zeta Function* defined as

$$\zeta(s) \equiv 1 + \frac{1}{2^s} + \frac{1}{3^s} + \cdots \equiv \sum_{1}^{\infty} \frac{1}{n^s}, \Re(s) > 1, s \equiv \alpha + i \cdot \beta$$

(*)

The domain of convergence of the *Riemann Zeta Function* can be extended for $\Re(s) < 1$ by the well-known process of analytical continuation.

The real zeros of the *Riemann Zeta-Function* occurring at $s = -2, -4, -6, \ldots$ are known as the trivial zeros. Other than these trivial zeros there are only complex zeros which are all to be found in the *critical strip*, $\Re(s) \in (0,1)$. The *RH* states that all non-trivial (e.g., complex) zeros of the *Riemann Zeta Function* lie on the *critical line*, $\Re(s) = \frac{1}{2}$. The *RH* is considered a pillar of the prime number theory, in that there is an important relationship between the zeros of the *Riemann Zeta-Function* and the distribution of the prime numbers. However, Riemann did not provide even a hint of a proof of his conjecture[1]. However, he proved that the number $N(T, \rho)$ of complex zeros with imaginary part $\tau$ between 0 and T (counted with multiplicities) is given by[2]

$$N(T, \rho) = \frac{T}{2\pi} \log \frac{T}{2\pi e} + \frac{7}{8} + S(T) + \mathcal{O}\left(\frac{1}{T}\right)$$

$$S(T) \equiv \frac{1}{\pi} ARG\left[\zeta\left(\frac{1}{2} + i \cdot T\right)\right], \zeta(\rho) = 0$$

(**)

Riemann also proved that $S(T) = \mathcal{O}(\log T)$ and suggested that the number $N_0(T, \rho_0)$ of zeros of $\zeta\left(\frac{1}{2} + i \cdot \tau\right), 0 < \tau \leq T$, seemed to be about,

$$\frac{T}{2\pi} \log \frac{T}{2\pi e},$$

(***)

which may have led him to his conjecture that all of zeros of the *Zeta-Function* in fact lie on the 1/2 vertical (*critical*) line.

In 1914 the distinguished British mathematician G. H. Hardy (1914) proved that an infinity of complex zeros of the *Zeta-Function* lie on the *critical line*. Subsequently Hardy and Littlewood (1918), Selberg (1942a) and (1942b), Levinson (1974) and Conrey (1989) have estimated the proportion of complex zeros of the *Riemann Zeta-Function* on the *critical line* to the number of complex zeros inside the *critical strip*,

---

[1] For a general background on the *RH*, see the survey article by Conrey (2003); also, the Clay Mathematics Institute website page https://claymath.org/millennium-problems/riemann-hypothesis has a good introduction. On the technical aspects of the *Riemann Zeta-Function*, see Bombieri (2010) for an excellent survey explaining in detail what is known about the *RH* problem and the many implications of a positive answer to the conjecture. The books by Edwards (1974) and Titchmarsh (1986) provide a comprehensive overview of the subject.
[2] Recall that $\mathcal{O}(f)$ denotes a quantity bounded in absolute value by a multiple of $f$.

$$N_0(T, \rho_0) \geq c_2, N(T, \rho), c_2 > 0, \zeta(\rho_0) = 0, \rho_0 \equiv \frac{1}{2} + i \cdot \tau$$

(#)

The (minimum) proportion $c_2$ of zeros on the *critical line* has been gradually improved to just over 40%.

Bohr and Landau (1914) have also shown that inside the arbitrarily thin strip defined by, $\Re(s) \in \left(\frac{1}{2} - \epsilon, \frac{1}{2} + \epsilon\right), \epsilon > 0$, lie "almost all" the non-trivial zeros of the *Riemann Zeta Function*. All in all, it is known that more than 40% of non-trivial zeros are simple and satisfy the *RH*. Van de Lune *et al.* (1986) showed that the first 1.5 billion zeros of the *Zeta-Function* - arranged by increasing positive imaginary part - are simple and satisfy the *RH*. With the help of high-speed computers more and more zeros of the *Zeta-Function* have been found to conform to the *RH*. Odlyzko (2001), who has computed the $10^{23}$-rd zero of the *Zeta-Function* and billions of its neighbours, notes that the fact that the first $10^{13}$ zeros are known to lie on the *critical line* (Gourdon, 2004) should not by itself be taken as sufficiently convincing evidence for *RH*[3]. In fact, the *RH* has so far been neither proved nor disproved. A list of important conjectures is a consequence of the *RH*, such as the *Lindelof Hypothesis* (see Edwards, 1974, ch. 9), the study of *L-Functions* (Selberg, 1991). Also, the *RH* has certain connections with Random Matrix Theory (RMT; Montgomery, 1973; Rudnick and Sarnak, 1996; Katz and Sarnak, 1999), resulting in a detailed model of the *Riemann Zeta-Function* and its value distribution; and with the inverse spectral problem for fractal strings (c.f. Lapidus and Maier, 1995). Some intriguing links with probability theory via the expectation in a moment of a Brownian bridge has been highlighted by Biane *et al.* (2001). The interplay between randomness and determinism is an important subject in the realm of prime numbers. Recently, LeClaire and Mussardo (2018) have shown that a random walk approach provides a key to establish the validity of the so-called *Generalised-Riemann-Hypothesis* (*GRH*) for the *Dirichlet L-Functions* of non-principal characters. *L-Functions* provide a generalisation of the *Riemann Zeta-Function* and the GRH conjecture states that the non-trivial zeros of all infinitely many *L-Functions* lie along the *critical line*, $\Re(s) = \frac{1}{2}$. According to LeClaire and Mussardo (2018) all *L-Functions* based on principal characters exactly share the same non-trivial zeros of the *Riemann Zeta-Function*. This means that the proof of the *RH* extends its reach to the *GRH* albeit limited to the *L-Functions* relative to principal characters[4].

---

[3] Some zeros around the $10^{21}$-st and $10^{22}$-nd are also known (Odlyzko, 1998). More recently, LeClaire (2017) has proposed a new algorithm for computing estimates of large zeros - as far as the $10^{100}$-th and even beyond with its associated error estimate which takes only few minutes to perform on a laptop using Mathematica.

In a fascinating recent essay Connes (2019) gives a through account of various approaches to the *RH* conjecture problem, and the state of efforts to pursue them, by navigating the many forms of the explicit *Riemann Zeta-Functions* and possible strategies to attack the problem, stressing the value of the elaboration of new concepts rather than "problem solving". In this work, however, we do take the more traditional "problem solving" approach without introducing new concepts.

---

[4] The *Riemann Zeta-Function* is a *Dirichlet L-Function* with a (trivial) principal character equal to one for every term of the infinite series.

In this study we present two proofs of the *RH* based on four important, well-known theorems concerning the properties of complex-valued analytic functions on a disk:

1) Jensen Theorem (formula) for the modulus of log-functions on a disk
2) Titchmarsh Theorem on the number of zeros in a disk
3) Rouché Theorem on the number of zeros on a disk
4) Riemann Mapping Theorem.

Applying the first three theorems to the zeros of the *Riemann Zeta-Function* requires a suitable mapping from the Unit Disk to the *critical strip*. Riemann Mapping Theorem provides the underpinning for the required analytic function. Perhaps one of the advantages of our approach is that it relies primarily – if not entirely - on a level of mathematical knowledge (advanced undergraduate) which makes it accessible to a relatively wider audience (familiar with complex numbers theory), and thereby allowing for a rapid check of its validity. Moreover, our approach does not use any functional equation form of the *RiemannZeta-Function*, other than leveraging its implied symmetry for non-trivial zeros on the *critical strip*.

Here is a short outline of the rest of this paper. Section 2 is devoted to reviewing the main properties of the *Riemann Zeta-Function* and its links with the *Eta-Function* integral form representation. Some important properties of the *Zeta-Function* integral representation are established. Section 3 briefly describes the connections of the non-trivial zeros of the *Zeta-Function* with those of its companion integral form representation. In Section 4, two proofs of the Riemann Hypothesis are presented; to focus on the key steps of the proofs, we had to relegate some of the lesser important technical details to an Appendix.

## 2. The *Riemann Zeta-Function* and the *Dirichlet Eta-Function* Integral Form

The *Riemann Zeta-Function*, denoted ζ(s), is introduced by means of the *Dirichlet Eta-Function* - also known as the alternating *Zeta-Function* (a type of *L-Function* with only a principal character of 1) - which can be expressed as follows (a special case of Dirichlet series):

$$\eta(s) \equiv (1 - 2^{1-s})\zeta(s) = \sum_{1}^{\infty} \frac{(-1)^{n+1}}{n^s}, \Re(s) > 0, s \in \mathbb{C}$$

(1)

where ζ(s) denotes the *Riemann Zeta Function* and η(s) has its domain of convergence for $\Re(s) > 0$. Notice that inside the *critical strip* where $\Re(s) \in (0,1)$ the factor $(1 - 2^{1-s})$ never vanishes and has no poles. Therefore, we can use the *Eta Function* as analytical continuation of the *Riemann Zeta Function* for $\Re(s) \in (0,1)$. The analytical continuation of ζ(s) for $\Re(s) < 0$ can be obtained using the remarkable *Riemann Functional Equation*[5],

$$\zeta(1-s) = \Gamma(s) \frac{2}{(2\pi)^s} \cos\left(\frac{s\pi}{2}\right) \zeta(s), \Re(s) \in (0,1), s \in \mathbb{C}$$

(2)

where $\Gamma(s)$ denotes the complex valued *Gamma Function*. Using *Riemann Functional Equation* (2) it can be showed that the non-trivial zeros of *ζ(s)* are located symmetrically with respect to the *critical line* $\Re(s) = 1/2$, inside the *critical strip*, $0 < \Re(s) < 1$. Hence, to prove the *RH* it is sufficient to show that $\zeta(s)$ never vanishes, for example, on the upper half of the *critical strip*, $1/2 < \Re(s) < 1$.

---

[5] See Titchmarsh (1986), p. 13.

Based on the Euler form it is easy to show that $\Gamma(s)$ do not vanish on the *critical strip*, $\Re(s) \in (0,1)$. Using formula (6.1.25, p. 256) in Abramowitz and Stegun, (1964), we get

$$|\Gamma(s)| \equiv |\Gamma(\alpha + i \cdot \beta)| = |\Gamma(\alpha)| \sqrt{\prod_{n=0}^{\infty} \frac{1}{1 + \frac{\beta^2}{(n+\alpha)^2}}} > 0,$$

$$\Gamma(\alpha) \equiv \int_0^\infty t^{\alpha-1} e^{-t}\, dt, \qquad s \equiv \alpha + \beta \cdot i, \alpha \in (0,1), \beta \in \mathbb{R}$$

(3)

The *Gamma-Function* $\Gamma(s)$ is a nowhere-vanishing meromorphic function with poles at the non-positive integers and no other poles. Hence, $|\Gamma(s)|$ and $\Gamma(\alpha)$ are both strictly positive and finite on the *critical strip*[6].

We refer to analytic continuation whenever using the *Riemann Zeta Function* for $\Re(s) \leq 1$. The *Dirichlet Eta-Function* (1) can be expressed in its integral form as follows[7]:

$$\Gamma(s)\eta(s) = \Gamma(s)(1 - 2^{1-s})\zeta(s) = \int_0^\infty \frac{x^{s-1}}{e^x + 1}\, dx, \Re(s) \in (0,1), s \in \mathbb{C}$$

(4)

---

[6] See Abramowitz and Stegun, (1964), formula 6.1.1, p. 255.
[7] See, for example, Sodlow (2005), formulas (17)-(18).

***Lemma 1***: The *Eta Function* (4) is well-defined as an improper integral and is finite if and only if,

$$F(s) \equiv \int_0^\infty \frac{x^{s-1}}{e^x + 1} dx, \Re(s) \in (0,1), s \in \mathbb{C}$$

(5)

is bounded, namely if its modulus is bounded by a finite value[8], denoted by *M*,

$$|F(s)| \equiv \left| \int_0^\infty \frac{x^{s-1}}{e^x + 1} dx \right| \leq M < +\infty, \Re(s) \in (0,1)$$

(6)

***Proof***: see Appendix, ***Proposition 1A*** where we compute the value of *M* as a function of $\Re(s)$ in closed form (see below, eq. 7).

∎

Moreover, we know from (3) that $|\Gamma(s)|$ is bounded on the *critical strip*. Therefore, the *Eta-Function* (4) is bounded as well. The upper bound of $|F(s)|$, *M*, is given by,

$$M(\alpha) \equiv \frac{1}{2\alpha} + e^{-1}, \alpha \equiv \Re(s), \alpha \in (0,1)$$

(7)

Hence, $M \in [1/2 + e^{-1}, +\infty)$. By restricting the domain to the upper half of the *critical strip*, $\Re(s) \in \left[\frac{1}{2}, 1\right]$, the supremum of $M(\alpha)$ is found at $\alpha = \frac{1}{2}$,

$$\left| \int_0^\infty \frac{x^{s-1}}{e^x + 1} dx \right| < 1 + e^{-1} \cong 1.36788, \Re(s) \in \left[\frac{1}{2}, 1\right]$$

(8a)

---

[8] See Apostol (1974), theorem 10.33, p. 276 on the existence of the improper integral.

since $M(\alpha)$ is a strictly decreasing function on the *critical strip*, as suggested by a direct inspection of (7). On the lower half of the *critical strip*, $\Re(s) \in \left(0, \frac{1}{2}\right)$, $M(\alpha)$ grows without bound ($M(\alpha) \to \infty$ for $\alpha \to 0$). The existence of a (finite) upper bound for (the modulus of) $F(s)$ on the (upper) half of the *critical strip* such that,

$$|F(s)| \leq F\left(\frac{1}{2}\right) < M\left(\frac{1}{2}\right)$$

(8b)

turns out to be a crucial insight in our proofs of the *Riemann Hypothesis* (see below). Thus, from now on we focus our attention on $|F(s)|$ solely on the upper half of the *critical strip*, $\Re(s) \in \left[\frac{1}{2}, 1\right]$. With (8b) we claim that we can improve upon the upper bound, $M\left(\frac{1}{2}\right) = 1 + e^{-1}$, on $\Re(s) \in \left[\frac{1}{2}, 1\right]$ and such improvement turns out to coincide with the $F(s)$ evaluated on the *critical line* at $s = \frac{1}{2}$,

$$F\left(\frac{1}{2}\right) \equiv \int_0^\infty \frac{x^{-\frac{1}{2}}}{e^x + 1} dx > 0$$

(8c)

**Remark 1**: The sharper bound (8c) is proved to be valid by inspecting the following chain of inequalities,

$$|F(s)| \equiv \left|\int_0^\infty \frac{x^{s-1}}{e^x + 1} dx\right| \leq \int_0^\infty \left|\frac{x^{s-1}}{e^x + 1}\right| dx \leq \int_0^\infty \frac{x^{\alpha-1}}{e^x + 1} dx, \alpha \equiv \Re(s), \alpha \in \left[\frac{1}{2}, 1\right]$$

(9)

and proving that the integral function

$$M^*(\alpha) \equiv \int_0^\infty \frac{x^{-\frac{1}{2}}}{e^x + 1} x^{\left(\alpha - \frac{1}{2}\right)} dx > 0, \alpha \in \left[\frac{1}{2}, 1\right]$$

(10a)

to be monotonically (strictly) decreasing in $\alpha$ and

$$M^*\left(\frac{1}{2}\right) = F\left(\frac{1}{2}\right) < M\left(\frac{1}{2}\right) = 1 + e^{-1}$$

(10b)

We start by showing that the first derivative of $M^*(\alpha)$ is negative,

$$\frac{dM^*}{d\alpha} = \int_0^\infty \frac{x^{-\frac{1}{2}}}{e^x + 1} x^{\left(\alpha - \frac{1}{2}\right)} \log(x)\, dx < 0, \alpha \in \left[\frac{1}{2}, 1\right]$$

(10c)

while its second derivative is positive

$$\frac{d^2 M^*}{d\alpha^2} = \int_0^\infty \frac{x^{-\frac{1}{2}}}{e^x + 1} x^{\left(\alpha - \frac{1}{2}\right)} [\log(x)]^2\, dx > 0, \alpha \in \left[\frac{1}{2}, 1\right]$$

(10d)

with the aid of the following,

<u>Lemma 2</u>: Let $g: \mathbb{R} \mapsto \mathbb{R}$ be a twice differentiable function then if $[a, b]$ is an interval on which $g''(x) > 0$ the function $g$ is convex on $[a, b]$, that is for $x < y \in [a, b]$ we have,

$$g[tx + (1 - t)y] \leq tg(x) + (1 - t)g(y), \quad t \in [0,1]$$

(11a)

Thus, informally speaking, chords between points on the graph of $g$ lie above the graph itself. Also, $g'$ is a monotonically increasing.

<u>Proof</u>: see Appendix, **Proposition 2A**.

∎

<u>Remark 2</u>: we can apply **Lemma 2** to $M^*(\alpha)$ on $\left[\frac{1}{2}, 1\right]$ since it is a twice differentiable function, as shown in (10c) and (10d), with positive second derivative – the latter assertion is also evident in that the integrand in (11b) is positive. Hence, we have that,

$$M^*(\alpha) \leq t \cdot M^*\left(\frac{1}{2}\right) + (1-t) \cdot M^*(1), \alpha = t \cdot \frac{1}{2} + (1-t) \cdot 1, \quad t \in [0,1]$$

(11b)

∎

Moreover, proving that

$$M^*\left(\frac{1}{2}\right) > M^*(1)$$

(12a)

would establish that,

$$M^*(\alpha) \leq M^*\left(\frac{1}{2}\right), \forall \alpha \in \left[\frac{1}{2}, 1\right]$$

(12b)

in that $M^*(\alpha)$ is monotonically decreasing. Also, notice that (8c) and (10a) imply that,

$$M^*\left(\frac{1}{2}\right) = F\left(\frac{1}{2}\right) \equiv \int_0^\infty \frac{x^{-\frac{1}{2}}}{e^x + 1} dx,$$

(12c)

To show that $M^*\left(\frac{1}{2}\right)$ is a tighter bound than $M\left(\frac{1}{2}\right) = 1 + e^{-1}$, and thereby confirming that (10b) holds, we compute the following integrals[9],

$$M^*\left(\frac{1}{2}\right) = 1.07215 < 1 + e^{-1} = 1.36788,$$

(13a)

and,

$$M^*(1) = \log(2) = 0.69315$$

(13b)

---

[9] Computed with the aid of Wolfram-Alpha (Mathematica) software (5 digits rounding).

In addition, $M^*\left(\frac{1}{2}\right)$ is an upper bound for the right-hand side of (12),

$$M^*(\alpha) \leq t \cdot M^*\left(\frac{1}{2}\right) + (1-t) \cdot M^*(1) \leq M^*\left(\frac{1}{2}\right), \alpha = t \cdot \frac{1}{2} + (1-t) \cdot 1, \quad t \in [0,1]$$

(14)

Being a convex function $M^*(\alpha)$ should have a monotonically increasing first derivative. Computing its values on the boundary of $\left[\frac{1}{2}, 1\right]$ we get,

$$\left.\frac{dM^*}{d\alpha}\right|_{\alpha=1/2} = \int_0^\infty \frac{x^{-\frac{1}{2}}}{e^x + 1} \log(x)\, dx = -1.76259$$

$$\left.\frac{dM^*}{d\alpha}\right|_{\alpha=1} = \int_0^\infty \frac{1}{e^x + 1} \log(x)\, dx = -\frac{1}{2}(\log 2)^2 = -0.240227$$

(15)

A negative value on the upper bound $\alpha = 1$ for $\frac{dM^*}{d\alpha}$ and the monotonically increasing property allows us to establish that the first derivative of $M^*(\alpha)$ must be negative on the interval $\left[\frac{1}{2}, 1\right]$, in that a change of sign for $\frac{dM^*}{d\alpha}$ would imply that there exists a subinterval in $\left[\frac{1}{2}, 1\right]$ where $\frac{dM^*}{d\alpha}$ should be decreasing rather than increasing. Such pattern would lead to a contradiction in that $\frac{dM^*}{d\alpha}$ is monotonically increasing having positive derivative (see eq. 10d). Furthermore, recalling (9), (12b) and (13a) we can assert that $F\left(\frac{1}{2}\right)$ provides the required upper bound,

$$|F(s)| \leq F\left(\frac{1}{2}\right) = M^*\left(\frac{1}{2}\right), \forall \Re(s) \in \left[\frac{1}{2}, 1\right]$$

(16)

**3. The Non-Trivial Zeros of the *Riemann Zeta-Function*.**

We can now proceed to make use of eqs. (1), (4) and (5) in writing the *Riemann Zeta-Function* in its integral form,

$$\zeta(s) = (1 - 2^{1-s})^{-1} \eta(s) = \frac{1}{(1 - 2^{1-s})\Gamma(s)} F(s), \quad \Re(s) \in (0,1)$$

(17a)

exploiting the *Eta-Function integral form*. Eq. (17a) is the *Master Equation* of this paper in building our proofs of the RH. As already stated, its denominator is non-vanishing on the *critical strip*,

$$(1 - 2^{1-s})\Gamma(s) \neq 0, \quad \Re(s) \in (0,1)$$

(17b)

Thus, we can ignore it in searching for the zeros of the *Riemann Zeta-Function,* while concentrating our attention on the integral function *F(s)* as defined in eq. (5),

$$F(s) \equiv \int_0^\infty \frac{x^{s-1}}{e^x + 1} dx = 0, \quad \Re(s) \in (0,1), s \equiv \alpha + i \cdot \beta$$

(17c)

Notice that *F(s)* can be interpreted as the Mellin transform[10] of the real valued function, $\frac{1}{e^x+1}, x \geq 0$. Locating the zeros of $\zeta(s)$ in the *critical strip* is equivalent to finding the zeros of the Mellin transform for $\frac{1}{e^x+1}, x \geq 0$. Furthermore, $F(s)$ is analytic in that it is the product of two analytic function, $\Gamma(s)$ and $\eta(s)$ as shown in eq. (4).

To recap our basic conclusion, $\zeta(s)$ and $F(s)$ share the same zeros in the *critical strip*. Thus, we can state the following corollaries (omitting their proofs),

***Corollary 1***:

$$\zeta(s) = 0 \Leftrightarrow F(s) \equiv \int_0^\infty \frac{x^{s-1}}{e^x + 1} dx = 0, \Re(s) \in (0,1)$$

---

[10] See Bleistein and Handelsman (1986), chapters 4-6, for an introduction to the theory and applications of Mellin transform.

(17d)

∎

and,

<u>*Corollary 2*</u>: by virtue of *Corollary 1* we can also assert that,

$$\zeta(s) \neq 0 \Leftrightarrow F(s) \neq 0$$

(17e)

∎

Thus, we focus heretofore on the zeros of *F(s)* in that they coincide with the non-trivial zeros of the *Riemann Zeta-Function*, $\zeta(s)$.

## *4. Two Proofs of the Riemann Hypothesis*

We state the *RH* as follows,

<u>*Theorem 1*</u> (*Riemann Hypothesis*): the non-trivial zeros of the *Riemann Zeta-Function* - e.g., the zeros inside the *critical strip* $\Re(s) \in (0,1)$ – lie on the vertical line of the complex plane, $\Re(s) = \frac{1}{2}$

$$\zeta(s) = 0 \Rightarrow \Re(s) = \frac{1}{2}, \Im(s) \neq 0, \forall \Re(s) \in (0,1)$$

(18)

where $\Im(s)$ denotes the imaginary part of *s*. We give two distinct proofs of the *RH*.

<u>*Proof N. 1*</u>: we need to recall two related theorems: Jensen's theorem and its celebrated formula and a theorem by Titchmarsh (1968) on the number of zeros of a complex function inside a disk. We start with Jensen's theorem formulation with ***Lemma 3*** (see Alhfors, 1953; Edwards, 1974, p. 40), let *f(z)* be a function which is defined and analytic throughout a disk $D(0,R) \equiv \{|z| \leq R,\ z \in \mathbb{C}, R > 0\}$. Suppose that *f(z)* has no zeros on the boundary circle $|z| = R$ and that inside the disk it has the zeros, $z_1, z_2, z_3, \ldots, z_n$ (where a zero of order *k* is included *k* times in the list). Suppose finally that $f(0) \neq 0$ then,

$$\log \left| f(0) \prod_{i=1}^{n} \frac{R}{z_i} \right| = \frac{1}{2\pi} \int_0^{2\pi} \log \left| f(Re^{i\theta}) \right| d\theta,$$

$$z \equiv \alpha_z + \beta_z \cdot i = Re^{i\theta}, R = |z|,$$

$$\alpha_z = R \cos \theta,\ \beta_z = R \sin \theta,$$

(19a)

where $\theta$ is the principal argument of *z* and *R* its modulus representing the polar coordinates of complex variable, *z*. If *f(z)* has no zeros inside the disk, then equation (19a) is merely,

$$\log |f(0)| = \frac{1}{2\pi} \int_0^{2\pi} \log \left| f(Re^{i\theta}) \right| d\theta$$

(19b)

∎

***Lemma 4*** (Titchmarsh's Number of Zeros Theorem; Titchmarsh, 1968); Conway, 1978, p. 282). Let $f(z), z \in \mathbb{C}$, be analytic in an open region containing $\overline{D}(0,R) \equiv \{z \in \mathbb{C}: |z| < R, R > 0\}$. Let $|f(z)| \leq M$ in the compact disk $D(0,R) \equiv \{z \in \mathbb{C}: |z| \leq R\}$ and suppose $f(0) \neq 0$. Then for $0 < \delta < 1$ the number of zeros of $f(z)$ on the disk $D(0, \delta R) \equiv \{z: |z| \leq \delta R\}$, denoted $N_0[f(z); z \in \overline{D}(0,R)]$, is less than or equal to

$$\frac{1}{\log(1/\delta)} \log \frac{M}{|f(0)|}, 0 < \delta < 1$$

(20a)

∎

**_Remark 3:_** notice that $f(z)$ has no zeros on the disk $D(0, \delta R)$,

$$N_0[f(z); \forall z \in D(0, \delta R)] = 0, \delta \in (0,1)$$

(20b)

if and only if expression (20a) counting the number of zeros is strictly less than one,

$$\frac{1}{\log(1/\delta)} \log \frac{M}{|f(0)|} < 1, \delta \in (0,1)$$

(20c)

which can indeed be simplified as,

$$\delta M < |f(0)|, \delta \in (0,1)$$

(20d)

Since _Lemma 4_ implies that

$$|f(0)| \leq M$$

(20e)

by virtue of (20d) and (20e) we must have,

$$\delta < \frac{|f(0)|}{M} \leq 1, \delta \in (0,1)$$

(20f)

Showing that (20d) holds is one of the key steps in our proof of the _RH_.

∎

In applying _Lemma 4_ – which is based on the well-known Jensen's theorem (_Lemma 3_) – it is convenient to perform the following change of variable for F(s),

$$\omega \equiv s - \frac{1}{2}, \Re(s) \in \left[\frac{1}{2}, 1\right], \Re(\omega) \in \left[0, \frac{1}{2}\right]$$

(21)

which, considering (with no loss of generality) only the upper half of the imaginary axis, $\Im(s) \geq 0$, yields

$$F_\omega(\omega) \equiv \int_0^\infty \frac{x^{-\frac{1}{2}}}{e^x + 1} x^\omega \, dx, \Re(\omega) \in \left[0, \frac{1}{2}\right], \Im(\omega) \geq 0,$$

(22)

In applying Titchmarsh's theorem (***Lemma 3***) to the zeros of $F_\omega(\omega)$ we follow closely Conway (1978; pp. 280-282) argument. However, Titchmarsh's theorem is valid if the domain of $f(z)$ is a Disk with a positive radius (say, equal to 1). Thus, to be able to apply ***Lemma 4***, we should make an appropriate change of variable for $F_\omega$ linking the open Unit Disk, $\bar{D}(0,1)$, to the vertical lower-half strip $\Re(\omega) \in \left(0, \frac{1}{2}\right)$, corresponding to the upper-half *critical strip*, $\Re(s) \in \left(\frac{1}{2}, 1\right)$. We build such bridge by introducing the following analytical function,

$$\phi: \bar{D}(0,1) \equiv \{z \in \mathbb{C}: |z| < 1\} \mapsto \bar{Z}_\phi \equiv \left\{\omega \in \mathbb{C}: \Re(\omega) \in \left(0, \frac{1}{2}\right), \Im(\omega) > 0\right\}$$

(24a)

Such analytic map exists as the Riemann Mapping Theorem assumptions are fulfilled for both open sets, $\bar{D}(0,1)$ and $\bar{Z}_\phi$ (see Appendix, Proposition 3A). We define the analytic function $\phi(z; b)$ as given by

$$\phi(z; b) \equiv \phi_\alpha(\alpha_z, \beta_z; b) + \phi_\beta(\alpha_z, \beta_z; b) \cdot i = \frac{1}{4} - \frac{1}{2\pi} i \cdot \text{Log}\left[\frac{1 + \theta(z; b)}{1 - \theta(z; b)}\right], z \in D(0,1)$$

(24b)

such that,

$$\phi_\alpha(\alpha_z, \beta_z; b) \equiv \frac{1}{4} + \frac{1}{2\pi} Arg\left[\frac{1 + \theta(z; b)}{1 - \theta(z; b)}\right]$$

$$\phi_\beta(\alpha_z, \beta_z; b) \equiv -\frac{1}{2\pi} Log\left|\frac{1 + \theta(z; b)}{1 - \theta(z; b)}\right|$$

(24c)

with,

$$\theta(z; b) \equiv \frac{z - b \cdot i}{1 + z \cdot b \cdot i}, 0 < b < 1$$

$$z \equiv \alpha_z + \beta_z \cdot i, \sqrt{\alpha_z^2 + \beta_z^2} < 1$$

(24d)

Notice that $\theta(z; b)$ is a conformal self-map of the open Unit Disk (see Gamelin, 2001, pp. 289-292 for details) indexed by parameter $b$. Importantly, for $z = 0$ we get,

$$\phi(0; b) = \phi_\alpha(0; b) = \frac{1}{4} + \frac{1}{2\pi} Arg\left[\frac{1 - b \cdot i}{1 + b \cdot i}\right], \phi_\beta(0; b) = 0, \forall b \in (0,1)$$

(25a)

with

$$\lim_{b \to 1} \phi(0; b) = \lim_{b \to 1} \phi_\alpha(0; b) = \frac{1}{4} + \frac{1}{2\pi} Arg\left[\frac{1 - 1 \cdot i}{1 + 1 \cdot i}\right] = \frac{1}{4} - \frac{1}{4} = 0$$

(25b)

Moreover, one can check the real part of $\phi$ maps on the *critical strip*,

$$\phi_\alpha(\alpha_z, \beta_z; b) \in \left(0, \frac{1}{2}\right), \forall b \in (0,1)$$

(26a)

in that it can be shown that

$$Arg\left[\frac{1 + \theta(z; b)}{1 - \theta(z; b)}\right] \in \left(-\frac{\pi}{2}, \frac{\pi}{2}\right), \forall b \in (0,1)$$

(26b)

***Remark 4***: parameter $b$ plays a key role in allowing $\phi(z,b)$ to fully cover the open (shifted upper) *critical strip* $\bar{Z}_\phi$, in that $\phi(0,b)$ can be made arbitrarily close to zero and thereby letting $F_\omega(0)$ arbitrarily close to $M^*\left(\frac{1}{2}\right)$ (see below, eq. 28a).

∎

Bearing in mind the change of variable $\phi(z;b)$, we can establish the following sequence of identities,

$$F_\omega(\omega) \equiv F_\omega[\phi(z;b)] \equiv F_z^\phi(z;b), \omega \equiv \phi(z;b) \in Z_\phi, z \in \bar{D}(0,1)$$

(27)

Moreover, we can make $F_z^\phi(0;b)$ arbitrarily close to the real value on the *critical line*,

$$\lim_{b \to 1} F_z^\phi(z;b) = F_\omega\left(\lim_{b \to 1} \phi_a(0;b)\right) = \lim_{\omega \to 0} F_\omega(\omega) = F\left(\frac{1}{2}\right) = M^*\left(\frac{1}{2}\right)$$

(28a)

By virtue of inequality (12) and the mapping on the *critical strip* obtained by $\phi(z;b)$, we can assert that $M^*\left(\frac{1}{2}\right)$ is the upper bound for $F_z^\phi$ on the Unit Disk,

$$\left|F_z^\phi(z;b)\right| = |F_\omega(\omega)| \leq M^*\left(\frac{1}{2}\right), \omega \in Z_\phi, z \in \bar{D}(0,1)$$

(28b)

Thus, all conditions are fulfilled to apply Titchmarsh's theorem (***Lemma 3***) to the analytic function $F_z^\phi(z;b)$ on the Unit Disk, $\bar{D}(0,1)$, e.g., letting $R = 1$. In counting the number of zeros of $F_z^\phi(z;b)$ on $\bar{D}(0,1)$, $N_0$, we need to make sure that inequality (20d) – and thereby (20f)– holds,

$$\delta \cdot M^*\left(\frac{1}{2}\right) < \left|F_z^\phi(0;b)\right| = F_z^\phi(0;b) = F_\omega(\phi_\alpha(0;b)), \delta \in (0,1), b \in (0,1)$$

(29a)

with

$$|F_z^\phi(0;b)| \leq M^*\left(\frac{1}{2}\right), \forall b \in (0,1)$$

(29b)

and thereby asserting that no zeros can be found in the *critical strip* other than those on the *critical line*. We posit the following real-valued function,

$$G(b) = F_z^\phi(0,b), \forall b \in (0,1)$$

$$G(b) \equiv F_\omega(\omega_0(b)), \omega_0(b) \equiv \phi_\alpha(0,b) = \frac{1}{4} + \frac{1}{2\pi}Arg\left[\frac{1-b\cdot i}{1+b\cdot i}\right]$$

(30a)

with

$$G(b) \equiv \int_0^\infty \frac{x^{-\frac{1}{2}}}{e^x+1}x^{\omega_0(b)}\,dx = M^*(\alpha), \alpha \in \left(\frac{1}{2},1\right), \omega_0(b) = \alpha - \frac{1}{2}, \quad \forall b \in (0,1)$$

(30b)

and

$$G(b) \leq F_\omega(0) = M^*\left(\frac{1}{2}\right), \lim_{b \to 1} G(b) = F_\omega(0), \forall b \in (0,1)$$

(30c)

We claim that $G(b)$ is monotonically increasing on $(0,1)$ in that,

$$\frac{dG(b)}{db} = G'(\omega_0)\frac{d\omega_0}{db} > 0, G'(\omega_0) < 0, \frac{d\omega_0}{db} < 0, \quad b \in (0,1)$$

(31)

as we know that $G(b)$ replicates by construction the value of $M^*(\alpha)$ on $\left(\frac{1}{2},1\right)$ (see eq. 30b), which is monotonically decreasing (see eq. 10c) as claimed in eq. (31). Since we can argue that

$$\frac{d\omega_0}{db} = \frac{1}{2\pi} \frac{dArg\left[\frac{1-b\cdot i}{1+b\cdot i}\right]}{db} = -\frac{i \cdot Delta\left(\frac{b+i}{i-b}\right)}{(b-i)^2} < 0, \forall b \in (0,1)$$

(32)

where $Delta(.)$ denotes the *Dirac Delta-Function*, we can assert that *G(b)* is monotonically increasing as stated in eq. (31).

Hence for any given $\delta \in (0,1)$ we can always pick a sufficiently large value, $b(\delta)$, close (or equal) to 1, such that,

$$\delta < \frac{G(b(\delta))}{M^*\left(\frac{1}{2}\right)} \leq 1, \quad b(\delta) \in (0,1), \forall \delta \in (0,1)$$

(33a)

in that *G(b)* is a continuous increasing function approaching $M^*\left(\frac{1}{2}\right)$ as *b* approaches 1,

$$\lim_{b \to 1} \frac{G(b)}{M^*\left(\frac{1}{2}\right)} = 1$$

(33b)

Hence, we can draw the conclusion that inequalities (29a) and (29b) always hold by setting parameter *b* arbitrarily close to 1 and thereby we can assert that $F_z^\phi(z;b)$ has no zeros on $\overline{D}(0,1)$ as appropriately gauged by the values of $b(\delta) \in (0,1)$

$$N_0\left[F_z^\phi(z;b); z \in \overline{D}(0,1)\right] = 0, b \in (0,1)$$

(34)

**_Remark 5_**: we can prove by contradiction that (34) must hold. Suppose there is a zero for $F_\omega$ in the interior of the lower half critical strip,

$$F_\omega(\widehat{\omega}) = 0, \Re(\widehat{\omega}) \in \left(0, \frac{1}{2}\right)$$

(34a)

we can solve the following equation for the value of $\hat{b}$ such that,

$$\phi_\alpha(0;\hat{b}) = \widehat{\omega_0}, \hat{b} \in (0,1), \widehat{\omega_0} \in (0, \Re(\hat{\omega}))$$

(34b)

which is

$$\hat{b} = \frac{1}{i}\left(\frac{1-\hat{\vartheta}}{1+\hat{\vartheta}}\right), \hat{\vartheta} \equiv ARG^{-1}\left[\left(2\cdot\widehat{\omega_0} - \frac{1}{2}\right)\frac{1}{\pi}\right]$$

(34c)

so that we can argue that $\phi(z;\hat{b})$ is a sufficient mapping covering $\hat{\omega}$ inside the critical strip and thereby fixing the corresponding value $\hat{z}$ inside the Unit Disk for a given value $\hat{b}$,

$$\phi(\hat{z};\hat{b}) = \hat{\omega}, \hat{z} \in D(0,1)$$

(34d)

yielding the zero, $\hat{z}$, for $F_z^\phi$. If it is the case that

$$\hat{\delta} \cdot M^*\left(\frac{1}{2}\right) < F_z^\phi(0;\hat{b}), \hat{\delta} = \lceil \hat{z} \rceil < 1, \hat{b} \in (0,1)$$

(34e)

then we immediately reach the conclusion that (34e) fulfils the required condition (29a), namely (20d) holds, with

$$f(z) \equiv F_z^\phi(z;\hat{b}), \cdot M \equiv, M^*\left(\frac{1}{2}\right), \delta \equiv \hat{\delta}, \hat{\delta} \in (0,1), \hat{b} \in (0,1)$$

(34f)

so that $F_z^\phi(z;\hat{b})$ cannot have a zero on the Disk $D(0,\hat{\delta})$. Hence, (34) holds true. If (34e) does not hold, we can just raise the value of $\hat{b}$, by lowering $\widehat{\omega_0}$ accordingly, so that $\phi(0;\hat{b})$ becomes sufficiently close to zero and thereby $F_z^\phi(0;\hat{b})$ close enough to $M^*\left(\frac{1}{2}\right)$ such that (34e) now holds in that $\hat{\delta} \in (0,1)$. We now check that we can always find a value $\hat{b} \in (0,1)$ consistent with inequality (34e). We proceed by considering the following Taylor series approximation around $b = 1$, arrested at the first order term, for $G(b)$,

$$G(b) = G(1) - \left.\frac{dG(b)}{db}\right|_{b=1} (1-b) + \mathcal{O}(b^2) =$$

$$=\cdot M^*\left(\frac{1}{2}\right) - \left[\left(\int_0^\infty \frac{x^{-\frac{1}{2}}}{e^x + 1} \log(x)\, dx\right)\frac{1}{2} Delta(-i)\right](1-b) + \mathcal{O}(b^2)$$

(34g)

with $\mathcal{O}(b^2)$ denoting the higher order terms.

We invert $\omega = \varphi[\theta(z;b)] \equiv \phi(z;b)$ (see, 24b-24c) and obtaining $\varphi^{-1}$ so that we can get $\theta$ to be a function of $\omega$. From (24d) we can find $z$ as function of $\theta$ and compute its modulus,

$$|z| \equiv H(\theta;b) = \left[\frac{b^2 + |\theta|^2 + 2b \cdot \mathfrak{Im}(\theta)}{1 + b^2|\theta|^2 + 2b \cdot \mathfrak{Im}(\theta)}\right]^{\frac{1}{2}}, \varphi^{-1}(\omega) = \theta, |\theta| \in (0,1), b \in (0,1)$$

(34h)

Let us consider the Taylor expansion for $H(\hat{\theta};b)$ as a function of $b$ by fixing the value $\hat{\theta} = \varphi^{-1}(\hat{\omega})$ for the given root $\hat{\omega}$,

$$\hat{\theta} = \varphi^{-1}(\hat{\omega}) = \frac{\psi(\hat{\omega}) - 1}{\psi(\hat{\omega}) + 1}, \psi(\hat{\omega}) \equiv \exp\left[\left(\hat{\omega} - \frac{1}{4}\right)\frac{2\pi}{i}\right]$$

$$|\hat{z}| \equiv H(\hat{\theta};b) = 1 - \left.\frac{dH(\hat{\theta};b)}{db}\right|_{b=1}(1-b) + \mathcal{O}(b^2) =$$

$$= 1 - \frac{2\left(1 - |\hat{\theta}|^2\right)}{1 + |\hat{\theta}|^2 + 2 \cdot \mathfrak{Im}(\hat{\theta})}(1-b) + \mathcal{O}(b^2), |\hat{z}| < 1$$

(34i)

Recalling (29a), for values of $b$ sufficiently close to 1 and fixing, $\hat{\delta} = [\hat{z}]$, so that $[\hat{z}] \in D(0, \hat{\delta})$, we can substitute the linear terms of the Taylor expansions (34g) and (34i) in (34e), with $G(b) = F_z^\phi(0; b)$, and check if the inequality holds,

$$\left[1 - \frac{2\left(1 - |\hat{\theta}|^2\right)}{1 + |\hat{\theta}|^2 + 2 \cdot \Im\mathfrak{M}(\hat{\theta})}(1-b)\right] \cdot M^*\left(\frac{1}{2}\right) <$$

$$< M^*\left(\frac{1}{2}\right) - \left[\left(\int_0^\infty \frac{x^{-\frac{1}{2}}}{e^x + 1}\log(x)\, dx\right)\frac{1}{2} Delta(-i)\right](1-b), \hat{\delta} \in (0,1), b \in (0,1)$$

(34j)

which yields,

$$\frac{2\left(1 - |\hat{\theta}|^2\right)}{1 + |\hat{\theta}|^2 + 2 \cdot \Im\mathfrak{M}(\hat{\theta})} \cdot M^*\left(\frac{1}{2}\right) > \left[\left(\int_0^\infty \frac{x^{-\frac{1}{2}}}{e^x + 1}\log(x)\, dx\right)\frac{1}{2} Delta(-i)\right] \cong$$

$$\cong -1.76259\frac{1}{2}4.66920, |\hat{\theta}| \in (0,1), M^*\left(\frac{1}{2}\right) > 0, |\hat{\theta}|^2 > \Im\mathfrak{M}(\hat{\theta})$$

(34k)

Since the denominator, $1 + |\hat{\theta}|^2 + 2 \cdot \Im\mathfrak{M}(\hat{\theta})$, is (strictly) positive as well as the numerator, $2\left(1 - |\hat{\theta}|^2\right)$, the fraction turns out to be positive and being $M^*\left(\frac{1}{2}\right)$ positive as well (e.g., eq. 13a), their product is positive whereas the term in the square bracket is negative, in that the integral is negative (e.g., eq. 15) and the Dirac Function value is positive[11]. Thus, we can assert that inequality (34k) holds, namely there exists a value $\hat{b}$ sufficiently close to 1 such that (34e) holds contradicting the assumption that there is a zero, $\hat{\omega}$, in the critical strip.

---

[11] $Delta(-i)$, the Dirac function evaluated at $-i$, is computed using Mathematica and rounded at the fifth digit.

∎

As there are no zeros on the lower half of the *critical strip* for $F_z(\omega)$, we can claim that $F(s)$ has no non-trivial zero on the upper half of the *critical line* and, by the symmetry of the functional equation (3), $F(s)$ is non-vanishing on the entire *critical strip* (other than the *critical line*). Thus, the *Riemann Zeta-Function* $\zeta(s)$ has no zeros as well. This completes our first proof of the *RH*.

∎

*Proof N. 2:* a version of the well-known Rouché's theorem is at the heart of our second proof. We recall the standard (symmetric) version of such theorem,

**Lemma 4 (asymmetric Rouche's Theorem[12])**: let the complex-valued function *f(s)* and *g(s)* be analytic in an open region *S* (a path-connected set), $S \subseteq \mathbb{C}$. Suppose that $\varrho$ is a simple closed curve in *S* enclosing a given region *D* of the complex plane. If

$$|g(s)| < |f(s)|, \forall s \in \varrho$$

(35)

then *f(s)* and *f(s) + g(s)* have the same number of zeros (counting multiplicities) inside the open subset $D \subset S$ bounded by $\varrho$.

∎

*Remark 6*: An equivalent result is the following (symmetric) version of Rouché's theorem. If

$$|f(s) + g(s)| < |f(s)| + |g(s)|, \forall s \in \varrho$$

---

[12] See Apostol (1974), p. 474-475, ex. 16-14. Rouché's theorem appear – sometimes under different guises – in many books giving an introduction on complex analysis (see, for example, Gamelin, 2001, p. 229).

(36)

then *f(s)* and *g(s)* have the same number of zeros (counting multiplicities) on *D*. In the literature the closed curve $\varrho$ goes under the name of Jordan curve. Notice that the boundary of the generic disk, $\varrho = \partial D(0, R)$, would be a special case of a simple closed curve enclosing $\bar{D}(0, R)$.

∎

Proofs of both versions of Rouché's theorem, as summarised by inequalities (35) and (36), are essentially the same (cf. Conway, 1978). More recently their validity has also been extended to arbitrary planar compacta taking the role of the open set *D* (*cf.* Narasimhan, 1985, p. 105).

A homotopic variant of Rouché's theorem on compact set has been recently presented by Mortini and Rupp (2014),

<u>*Lemma 5*</u> *(Rouche's Theorem for homotopic maps; Mortini and Rupp, 2014, th. 15, p. 7)*, let $K \subseteq \mathbb{C}$ be a compact set and $f(s), g(s) \in \mathcal{A}(K)$ where $\mathcal{A}(K)$ denotes the space of all functions continuous on K and holomorphic in the interior of K. Suppose that $f(s)$ and $g(s)$ are homotopic in $\mathcal{C}(\partial K, \mathbb{C}^*), \mathbb{C}^* \equiv \mathbb{C}\setminus\{0\}$ then $f(s)$ and $g(s)$ have the number of zeros on K,

$$N_0[f(s); K] = N_0[g(s); K], s \in K$$

(37)

∎

<u>*Remark 7*</u>: the homotopic version of Rouché's theorem yields another proof of the classic Rouché's theorem itself. We just need to note that the condition,

$$|f(s) + g(s)| < |f(s)| + |g(s)|, \forall s \in \partial K$$

(38)

implies that the function $L: [0,1] \times ðK \to \mathbb{C}^*$ given by

$$L(t,s) \equiv (1-t)f(s) - tg(s)$$

(39)

is a homotopy that connects $f(s)$ with $-g(s)$ inside $\mathcal{C}(ðK, \mathbb{C}^*)$.

∎

Since we will use the homotopy version of Rouché's theorem given in **Remark 7** we can skip the map $\phi(z)$ and proceed by directly considering the lower half of the *critical strip* on the upper half- plane (non-negative imaginary axis),

$$LHCS \equiv \left\{ \omega \in \mathbb{C} : \Re(\omega) \in \left[0, \frac{1}{2}\right], \Im(\omega) \geq 0 \right\}$$

(40)

We define a compact set inside *LHCS* by crossing the imaginary axis with a horizontal line setting at level, $\tau$, the height of a rectangle inside the *critical strip*,

$$K(\tau) \equiv \{\omega \in LHCS : 0 \leq \Im(\omega) \leq \tau, \quad \tau > 0\}$$

(41)

In applying the homotopic version of Rouché's theorem, we posit the following definition,

$$\begin{cases} f(\omega) \equiv F_\omega(\omega) \cdot L(\omega) \\ g(\omega) \equiv \lambda \cdot (\epsilon + \omega) \end{cases}, \quad \varepsilon, \lambda > 0, \forall \omega \in LHCS,$$

(42a)

with

$$L(\omega) \equiv \prod_{j=1}^{j=\rho(\tau)} \left[ \frac{\overline{\omega} + i \cdot \beta_j}{\omega - i \cdot \beta_j} \right], |L(\omega)| = 1, \forall \omega \in LHCS$$

(42b)

and $\omega = i \cdot \beta_j$ are the zeros of the Riemann Zeta Function on the (shifted) half *critical line* up to the height of the imaginary axis indexed by $\tau$,

$$F_\omega(i \cdot \beta_j) = 0, \beta_j > 0, \forall j = 1, \rho(\tau)$$

(42c)

with $\rho(\tau)$ denoting the number of zeros on the *critical line* up to the height, $\tau$. We argue that $L(\omega)$ has a unit modulus, $|L(\omega)| = 1, \forall \omega \in LHCS$, since the numerator, $(\bar{\omega} + i \cdot \beta_j)$, and denominator, $(\omega - i \cdot \beta_j)$, in each fraction of the product are complex conjugate of each other and therefore they have the same modulus. The factor $L(\omega)$ in the product $F_\omega(\omega) \cdot L(\omega)$ plays the role of a "neutraliser" of the zeros located on the *critical line*, guaranteeing that $f(\omega)$ is non-vanishing on that portion of the (truncated) boundary of the *critical strip*, $\partial K(\tau)$. Being that the case, we can argue more broadly that $f(\omega)$ is non-vanishing on the entire boundary $\partial K(\tau)$,

$$|f(\omega)| = |F_\omega(\omega) \cdot L(\omega)| > 0, \forall \omega \in \partial K(\tau)$$

(43)

in that $F_\omega(\omega)$ has no zeros on the whole interval of the real axis, $\forall \omega \in \left[0, \frac{1}{2}\right]$, as well as on the entire vertical axis, $\Re(\omega) = \frac{1}{2}$ which corresponds to the vertical line $\Re(s) \equiv \alpha = 1$. Moreover, we can rule out the *de-facto* non generic case with zeros located on the horizontal line $\Im(\omega) = \tau$,

$$F_\omega(\alpha_\omega + i \cdot \tau) = 0, \alpha_\omega \in \left[0, \frac{1}{2}\right]$$

(44)

in that the zeros of analytic functions are isolated points[13]. If, by chance, there were zeros with imaginary part equal to, say $\hat{\tau} > 0$, one can shift upward the horizontal line crossing the imaginary axis. There is always a neighbourhood of $\hat{\tau}$ – say a Disk with a positive radius - such that for an arbitrary $\varepsilon > 0$, $F_\omega(\alpha_\omega + i \cdot (\hat{\tau} + \varepsilon))$ is non-vanishing[14]. Moreover, we can argue that $g(\omega)$ is also a non-vanishing function,

$$|g(\omega)| \equiv \lambda \cdot \sqrt{(\alpha_\omega + \epsilon)^2 + \beta_\omega^2} > 0, \omega \equiv \alpha_\omega + i \cdot \beta_\omega, \alpha_\omega \geq 0, \epsilon > 0, \omega \in \partial K(\tau)$$

(45)

in that its real part, $(\alpha_\omega + \epsilon)$, is strictly positive.

Since both f($\omega$) and g($\omega$) are evidently analytic - f($\omega$) is a product of analytic functions and g($\omega$) is a linear affine transformation and do not vanish on the boundary, $\partial K(\tau)$, of a compact set, they fulfil the basic Rouché's theorem requirements. However, we still need to check that the (strict) inequality (38) holds true. By virtue of the triangle inequality theorem for complex variables we can assert that,

$$|F_\omega(\omega) \cdot L(\omega) + \lambda \cdot (\epsilon + \omega)| \leq |F_\omega(\omega)| + \lambda \cdot \sqrt{(\alpha_\omega + \epsilon)^2 + \beta_\omega^2}, \forall \omega \in \partial K(\tau)$$

(46a)

recalling that $|L(\omega)| = 1$. However, inequality (46a) is not sufficient *per se* as we need strict inequality for Rouché's theorem to be applicable, namely one should prove that,

$$|F_\omega(\omega) \cdot L(\omega) + \lambda \cdot (\epsilon + \omega)| < |F_\omega(\omega)| + \lambda \cdot \sqrt{(\alpha_\omega + \epsilon)^2 + \beta_\omega^2}, \forall \omega \in \partial K(\tau)$$

(46b)

---

[13] See Agarwal *et al.* (2011), Lecture 26, pp.177-182.
[14] If there were more than one value of $\alpha_\omega$ fulfilling eq. (44) one could always increase $\hat{\tau}$ by an arbitrary amount such that $F(\alpha_\omega + i \cdot \tau)$ is non-vanishing on the Disk centred in the point, $\alpha_\omega + i \cdot \hat{\tau}$, having the smallest radius.

To rule out the case of equality in (46a) and therefore asserting that

$$|F_\omega(\omega) \cdot L(\omega) + \lambda \cdot (\epsilon + \omega)| \neq |F_\omega(\omega)| + \lambda \cdot \sqrt{(\alpha_\omega + \epsilon)^2 + \beta_\omega^2}, \forall \omega \in \partial K(\tau)$$

(47a)

the following necessary and sufficient condition must be fulfilled,

$$F_\omega(\omega) \cdot L(\omega) \neq \vartheta(\omega) \cdot \lambda \cdot (\epsilon + \omega), \quad \vartheta(\omega) \neq 0, \forall \omega \in \partial K(\tau),$$

(47b)

for $\vartheta(\omega)$ being a real-valued function. In fact, we show that if (and only if)

$$F_\omega(\omega) \cdot L(\omega) = \vartheta(\omega) \cdot \lambda \cdot (\epsilon + \omega), \quad \vartheta(\omega) \neq 0, \omega \in \partial K(\tau)$$

(48a)

holds for some $\omega \in \partial K(\tau)$, then equality,

$$|F_\omega(\omega) \cdot L(\omega) + \lambda \cdot (\epsilon + \omega)| = |F_\omega(\omega)| + \lambda \cdot \sqrt{(\alpha_\omega + \epsilon)^2 + \beta_\omega^2}, \quad \omega \in \partial K(\tau)$$

(48b)

would be true[15].

We prove by contradiction that an appropriate choice of parameter $\lambda$ would prevent equality (48b) to holding on the boundary $\partial K(\tau)$. and thereby being able to argue that strict inequality (46b) is a valid statement. Let us suppose that for some $\widetilde{\omega} \in \partial K(\tau)$, there exists a value $\tilde{\vartheta}_\tau \neq 0$ such that (48b) holds, e.g.

$$F_\omega(\widetilde{\omega}) \cdot L(\widetilde{\omega}) = \tilde{\vartheta}_\tau \cdot \lambda \cdot (\epsilon + \widetilde{\omega}), \widetilde{\omega} \in \partial K(\tau),$$

(49a)

As a result, by taking the modulus of both side in (49a) we must have,

---

[15] See Appendix, **Proposition 4A**, for a proof that equality (48a) is a necessary and sufficient condition for equality (48b) to hold.

$$|F_\omega(\widetilde{\omega})| = |\tilde{\vartheta}_\tau| \cdot \lambda \cdot \sqrt{(\tilde{\alpha}_\omega + \epsilon)^2 + \tilde{\beta}_\omega^2}, |\tilde{\vartheta}_\tau| > 0, \widetilde{\omega} \in \partial K(\tau)$$

(49b)

with

$$|L(\widetilde{\omega})| = 1, |(\epsilon + \widetilde{\omega})| \equiv \sqrt{(\tilde{\alpha}_\omega + \epsilon)^2 + \tilde{\beta}_\omega^2}, \widetilde{\omega} \in \partial K(\tau)$$

(49c)

Without loss of generality, let us assume that $|\tilde{\vartheta}_\tau|$ is the lowest (positive) value (infimum) for $|\tilde{\vartheta}_\tau|$ on $\widetilde{\omega} \in \partial K(\tau)$ such that eq. (48) is supposed to hold. If we set parameter $\lambda_\tau$ such that

$$|\tilde{\vartheta}_\tau| \cdot \lambda_\tau \cdot \sqrt{(\tilde{\alpha}_\omega + \epsilon)^2 + \tilde{\beta}_\omega^2} > M^*\left(\frac{1}{2}\right) \Rightarrow \lambda_\tau > M^*\left(\frac{1}{2}\right) \Big/ \left[|\tilde{\vartheta}_\tau| \cdot \sqrt{(\tilde{\alpha}_\omega + \epsilon)^2 + \tilde{\beta}_\omega^2}\right] > 0,$$

(50a)

and substitute in the right-hand side of eq. (49b) the following value for parameter $\lambda_\tau$,

$$\tilde{\lambda}_\tau = \left[M^*\left(\frac{1}{2}\right) + \nu\right] \Big/ (|\tilde{\vartheta}_\tau| \cdot \epsilon), \nu > 0$$

(50b)

we would end up with a contradiction,

$$|F_\omega(\widetilde{\omega})| = \left[M^*\left(\frac{1}{2}\right) + \nu\right] \frac{\sqrt{(\tilde{\alpha}_\omega + \epsilon)^2 + \tilde{\beta}_\omega^2}}{\epsilon} > M^*\left(\frac{1}{2}\right), \widetilde{\omega} \in \partial K(\tau) \subset LHCS$$

(50c)

since we know that $M^*\left(\frac{1}{2}\right)$ is an upper bound for the modulus of $F_\omega(\omega)$ on the entire (shifted) lower half *critical strip*, *LHCS*. Hence, we can conclude that strict inequality (46b) always holds, in that (49a) can never be true if, as shown in (50b), the right choice of parameter, $\tilde{\lambda}_\tau$, is performed.

Thus, we can apply Rouché's theorem asserting that the number of zeros for $f(\omega)$ and $g(\omega)$ should coincide on $K(\tau)$,

$$N_0[f(\omega)] = N_0[g(\omega)], \omega \in K(\tau),$$

(51a)

recalling that on the boundary $\partial K(\tau)$ both functions are non-vanishing. Since it is evident that $g(\omega)$ is non-vanishing on $K(\tau)$, we should have,

$$N_0[g(\omega)] \equiv N_0[\tilde{\lambda} \cdot (\epsilon + \alpha_\omega) + i \cdot \tilde{\lambda} \cdot \beta_\omega] = 0, \alpha_\omega \geq 0, \tilde{\lambda} > 0, \epsilon > 0, \omega \in K(\tau)$$

(51b)

Taken together (51a) and (51b) imply that $f(\omega)$ is non-vanishing on $K(\tau)$ as well,

$$N_0[f(\omega)] \equiv N_0[F_\omega(\omega) \cdot L(\omega)] = 0, \omega \in K(\tau)$$

(52a)

thereby yielding the number of zeros for $F_\omega(\omega)$ in the open (truncated) *critical strip*,

$$N_0[F_\omega(\omega)] = 0, \omega \in \{K(\tau) \backslash \partial K(\tau)\}$$

(52b)

Assertion (52b) is proved by contradiction. Since on $\{K(\tau) \backslash \partial K(\tau)\}, L(\omega) \neq 0, f(\omega) \equiv F_\omega(\omega) \cdot L(\omega)$ can vanish if and only if $F_\omega(\omega) = 0$. But (52a) rules out the existence of a zero for $f(\omega)$ on $K(\tau)$ and *a fortiori* on $\{K(\tau) \backslash \partial K(\tau)\}$. Thus, asserting that $F_\omega(\omega)$ has a zero on $\{K(\tau) \backslash \partial K(\tau)\}$ would contradict (52a). Hence, (52b) need to hold, namely $F_\omega(\omega)$ is non-vanishing inside the (truncated) half *critical strip*, $\{K(\tau) \backslash \partial K(\tau)\}$.

We now extend the same assertion on the absence of zeros for $F_\omega(\omega)$ to the entire (open) *critical strip*, $\overline{LHCS}$. We proceed *by a reductio ad absurdum*, and so let us assume that there exist a zero for $F_\omega$, $\hat{\omega}$, which is not in the (truncated) open half *critical strip* $\{K(\tau) \backslash \partial K(\tau)\}$, but would still be in $\overline{LHCS}$,

$$F_\omega(\widehat{\omega}) = 0, \widehat{\omega} \equiv \widehat{\alpha}_\omega + i \cdot \widehat{\beta}_\omega, \widehat{\beta}_\omega > \tau, \widehat{\alpha}_\omega \in \left(0, \frac{1}{2}\right), \widehat{\omega} \in \overline{LHCS}$$

(53a)

By increasing the height of the (truncated) open half *critical strip* to just above the imaginary part of the zero, $\widehat{\beta}_\omega$

$$\widehat{\tau} = \widehat{\beta}_\omega + \varepsilon, \varepsilon > 0$$

(53b)

we can argue that (by construction) the rectangle with height $\widehat{\tau}$ includes the zero, $\widehat{\omega} \in \{K(\widehat{\tau}) \setminus \partial K(\widehat{\tau})\}$. However, $F_\omega(\widehat{\omega}) = 0, \widehat{\omega} \in \{K(\widehat{\tau}) \setminus \partial K(\widehat{\tau})\}$ would now contradict (52b) which holds for all positive $\tau$, including $\widehat{\tau}$. Hence, $F_\omega(\omega)$ is non-vanishing on the entire open (lower) half *critical strip* and thereby $F(s)$ has no non-trivial zeros on the open (upper) half *critical strip*. By the symmetry of the functional equation (3), $F(s)$ turns out to be non-vanishing on the lower half of the *critical strip* as well. Therefore, we can assert that $\zeta(s)$ cannot have any zero on the *critical strip* other than the *critical line* as well. Thus, the *RH* is again proved to be true.

.

## Acknowledgements

I would like to thank Michela Ablondi, Jacopo Montali, Andrew Tilling, Elisabetta Villani and Caterina Violi, for lively comments and interesting discussions, and a lot of moral support during the time of this work. Long walks with my grand-daughter Anna Isabel Tilling Violi on the Tuscan hills gave me the opportunity to free up my mind and concentrating on the important issues at a crucial stage of this work.

# Appendix

**_Proposition 1A_**: We prove that inequality (6) holds (**_Lemma 1_**).

**_Proof_**: The following chain of inequality can be checked by inspection to hold,

$$\left|\int_0^\infty \frac{x^{s-1}}{e^x+1} dx\right| \leq \int_0^\infty \left|\frac{x^{s-1}}{e^x+1}\right| dx \leq \int_0^\infty \frac{x^{\Re(s)-1}}{e^x+1} dx <$$

$$< \int_0^1 \frac{x^{\Re(s)-1}}{e^x+1} dx + \int_1^\infty \frac{1}{e^x+1} dx < \int_0^1 \frac{x^{\Re(s)-1}}{e^x+1} dx + \int_1^\infty e^{-x} dx$$

(1A)

We solve the first term in the last line by using integration by part where $\int v du = [uv] - \int u dv$ with,

$$du \equiv x^{\Re(s)-1} dx, v \equiv \frac{1}{e^x+1}$$

(2A)

which then yields,

$$\int_0^1 \frac{x^{\Re(s)-1}}{e^x+1} dx = \left[\frac{x^{\Re(s)}}{\Re(s)} \frac{1}{e^x+1}\right]_0^1 + \int_0^1 \frac{x^{\Re(s)} e^x}{\Re(s)(e^x+1)^2} dx <$$

$$< \left[\frac{x^{\Re(s)}}{\Re(s)} \frac{1}{e^x+1}\right]_0^1 + \int_0^1 \frac{e^x}{\Re(s)(e^x+1)^2} dx <$$

$$< \frac{1}{\Re(s)} \frac{1}{e^1+1} + \frac{1}{\Re(s)} \left(\frac{1}{2} - \frac{1}{e^1+1}\right)$$

(3A)

Therefore, by inserting (9) into (7) we get,

$$\left|\int_0^\infty \frac{x^{s-1}}{e^x+1} dx\right| < \frac{1}{\Re(s)} \frac{1}{e^1+1} + \frac{1}{\Re(s)} \left(\frac{1}{2} - \frac{1}{e^1+1}\right) + e^{-1} < +\infty, \Re(s) \in (0,1)$$

(4A)

and by setting,

$$M \equiv \frac{1}{\Re(s)} \frac{1}{e^1+1} + \frac{1}{\Re(s)}\left(\frac{1}{2} - \frac{1}{e^1+1}\right) + e^{-1} = \frac{1}{\Re(s)}\frac{1}{2} + e^{-1} < +\infty, \Re(s) \in (0,1)$$

(5A)

with $M > 0$. Moreover, if we restrict the domain to the upper half of the *critical strip*, the greatest upper bound $M$ can be fixed at

$$\left|\int_0^\infty \frac{x^{s-1}}{e^x+1} dx\right| < 1 + e^{-1} \cong 1.36788, \Re(s) \in \left[\frac{1}{2}, 1\right]$$

(6A)

Thus, we have proved in (4A) that inequality (6) holds and therefore the integral function $F(s)$ is bounded by a positive, finite constant $M$ in the upper half of the critical strip, as obtained in (6A).

∎

***Proposition 2A***: we prove the claims in ***Lemma 2*** asserting that $g$ is a convex function,

$$g[tx + (1-t)y] \leq tg(x) + (1-t)g(y), \forall x, y \in [a, b], t \in [0,1]$$

(7A)

and its first derivative, $g'$, is monotonically increasing, namely if $x < y$ we have that $g'(x) \leq g'(y)$. Moreover, if $g'(y) < 0$, we get,

$$g'[tx + (1-t)y] < 0, \forall\ t \in [0,1], \quad x < y, \quad x, y \in [a, b]$$

(8A)

namely g' is negative in the subinterval $[x, y]$. Furthermore, the inequality

$$tg(x) + (1-t)g(y) \leq \text{Max}[g(x), g(y)], \forall\ t \in [0,1]$$

evidently holds implying that,

$$g[tx + (1-t)y] \leq \text{Max}[g(x), g(y)]$$

(9A)

***Proof***: let $\xi \equiv tx + (1-t)y, t \in (0,1)$ be a point in the interval between *x* and *y*. Now the slope of the chord between (*x*, *g(x)*) and ($\xi$, *g($\xi$)*) is, by the Mean-Value Theorem, equal to *g'(s₁)* where *s₁* lies between *x* and $\xi$, while the slope of the chord between ($\xi$, *g($\xi$)*) and (*y*, *g(y)*) is equal to *g'(s₂)* for *s₂* between $\xi$ and *y*. Let us assume that *g* is convex, then it should be the case that,

$$g(\xi) < tg(x) + (1-t)g(y), t \in (0,1)$$

(10A)

which implies, by the Mean-Value Theorem, that the slopes of the chord have opposite sign as follows,

$$g'(s_2) > 0, g'(s_1) < 0$$

(11A)

However, if it were the case that *g* is not to be convex then (10A) does not hold and therefore it should be true that,

$$g(\xi) \geq tg(x) + (1-t)g(y), t \in (0,1)$$

(12A)

which would reverse the signs in (11A) yielding,

$$g'(s_2) \leq 0, g'(s_1) \geq 0$$

(13A)

which would lead one to assert that,

$$g'(s_1) - g'(s_2) > 0$$

(14A)

However, (14A) would lead to a contradiction (see below). By the Mean-Value Theorem for *g'(s)* applied to the values *s₁* and *s₂* it follows that there be an $s \in (s_1, s_2)$ such that,

$$g''(s) = (g'(s_2) - g'(s_1)) / (s_2 - s_1) > 0,$$

(15A)

in that $g''(s)$ is positive by assumption. As a result, the numerator in (15A) should be positive,

$$g'(s_2) - g'(s_1) > 0,$$

(16A)

since the denominator in (15A) is positive being the case that.

$$s_2 > s_1$$

(17A)

Thus, (16A) contradicts (14A) and thereby (12A). Hence, (10A) holds implying that $g$ is convex.

∎

***Proposition 3A***:(*Riemann Mapping Theorem; see Gamelin, 2001, p. 311-313*): Let $U$ be a simply connected open subset of $\mathbb{C}$ that is not all $\mathbb{C}$. Then $U$ is complex diffeomorphic to $\overline{D}(0,1)$.

∎

*Corollary 3A.1*: Let $U$ be a simply connected region and let $w \in U$, then, there is a unique bijective conformal transformation $f$ of $U$ onto the (open) Unit Disc, $\overline{D}(0,1)$, such that $f'(w) = 0, f'(w) > 0$.

∎

*Corollary 3A.2*: Any two simply connected domains have a conformal map between them. If $U$ is a convex set, then it is simply connected, so Riemann's mapping theorem ensure that there is a conformal map sending $U$ to $\overline{D}(0,1)$ and $\overline{D}(0,1)$ to $U$.

∎

*Remark 3A*: of course, knowing that a conformal map between two domains exists still leaves the challenge of constructing a right one. Such challenge cannot be (easily) avoided if there are constraints needed to be fulfilled on the adopted conformal map.

∎

*Proposition 4A*: we examine the relationship between two complex numbers, w and v, with respect to their modulus, $|w + v|, |w|$ and $|v|$. Triangular inequality implies $|w + v|$ is never larger than the sum $|w| + |v|$. To prove it, we write the triangle inequality,

$$|w| + |v| \geq |w + v|$$

(18A)

in terms of $w \equiv \alpha_w + i \cdot \beta_w$ and $v \equiv \alpha_v + i \cdot \beta_v$,

$$\sqrt{(\alpha_w)^2 + (\beta_w)^2} + \sqrt{(\alpha_v)^2 + (\beta_v)^2} \geq \sqrt{(\alpha_w + \alpha_v)^2 + (\beta_w + \beta_v)^2}$$

(19A)

squaring both sides and cancelling common terms leaves

$$2 \cdot \sqrt{(\alpha_w)^2 + (\beta_w)^2} \cdot \sqrt{(\alpha_v)^2 + (\beta_v)^2} \geq 2 \cdot \alpha_w \alpha_v + 2 \cdot \beta_w \beta_v$$

(20A)

Cancelling the 2's and squaring again, we get

$$[(\alpha_w)^2 + (\beta_w)^2][(\alpha_v)^2 + (\beta_v)^2] \geq (\alpha_w \alpha_v)^2 + (\beta_w \beta_v)^2 + 2\alpha_w \alpha_v \beta_w \beta_v$$

(21A)

or

$$(\alpha_w)^2 (\beta_v)^2 - 2\alpha_w \alpha_v \beta_w \beta_v + (\alpha_v)^2 (\beta_w)^2 \geq 0$$

(22A)

recognizing the left as a perfect square, we can write (22A) as

$$(\alpha_w \beta_v - \alpha_v \beta_w)^2 \geq 0$$

(23A)

which is clearly always true. Furthermore, we can find when equality,

$$|w| + |v| = |w + v|$$

(24A)

holds by noting that these can only be equal when,

$$(\alpha_w \beta_v - \alpha_v \beta_w)^2 = 0$$

(25A)

Clearly, this occurs if and only if,

$$\alpha_w \beta_v = \alpha_v \beta_w$$

(26A)

or

$$\frac{\alpha_w}{\alpha_v} = \frac{\beta_w}{\beta_v}$$

(27A)

If we let this common ratio be $\vartheta$, we find that our equality condition is,

$$\alpha_w = \vartheta \cdot \alpha_v \text{ and } \beta_w = \vartheta \cdot \beta_v,$$

(28A)

so, we must have

$$w = \vartheta \cdot v$$

(29A)

and therefore one (complex) number is just a real multiple, $\vartheta$, of the other.

∎